\UseRawInputEncoding

\documentclass[12pt]{amsart}
\usepackage{amsfonts,amsmath,amsthm,amssymb,indentfirst,url}

\usepackage[english,russian]{babel}

\renewcommand{\le}{\leqslant}
\renewcommand{\ge}{\geqslant}
\renewcommand{\leq}{\leqslant}

\newtheorem{theorem}{Теорема}
    \newtheorem{lemma}[theorem]{Лемма}
    \newtheorem{problem}[theorem]{Задача}
    \newtheorem{remark}[theorem]{Замечание}

\begin{document}

%Simple proofs of estimations of Ramsey numbers and of discrepancy
%A. Buchaev and A. Skopenkov
%In this expository note we present simple proofs of the lower bound of Ramsey numbers (Erd\"os theorem),
%and of the estimation of discrepancy.
%Neither statements nor proofs require any knowledge beyond high-school curriculum (except a minor detail).
%Thus they are accessible to non-specialists, in particular, to students.
%Our exposition is simpler than the standard exposition because no probabilistic language is used.
%In order to prove the existence of a `good' object we prove that the number of `bad' objects
%is smaller than the number of all objects.
%05-01; 05D10, 05D40

\title{Простые доказательства  \\  оценок чисел Рамсея и уклонения}
\author{А. Бучаев и А. Скопенков}
\thanks{А. Бучаев:  Московский Физико-Технический Институт.
\newline
А. Скопенков:
\texttt{https://users.mccme.ru/skopenko}.
Московский Физико-Технический Институт, Независимый Московский Университет.
\newline
Эта заметка возникла в ходе обсуждений в рамках семинаров по дискретному анализу
по курсу А.М. Райгородского на ФИВТ МФТИ.
Благодарим за полезные обсуждения Д.А. Колупаева, Г.М. Кучерявого, А.А. Печенкина, А.М. Райгородского
и анонимного рецензента.}
\date{}
\maketitle

\bigskip
Мы приводим простые доказательства теоремы  \ref{t:erd} Эрдеша о нижней оценке чисел Рамсея и теоремы \ref{t:rad} об оценке уклонения (см. формулировки  ниже).
Для понимания формулировок и доказательств не требуется знаний, выходящих за пределы школьной программы (кроме неравенства (*) в конце доказательства леммы \ref{l:one}, для которого требуется разложение экспоненты в ряд).
По сути наше изложение аналогично \cite[\S 1.1]{AS}, \cite[\S 3.2, \S 4.2]{R3}.
Однако оно проще для восприятия, поскольку не использует ненужного здесь вероятностного языка
(см. подробнее замечание \ref{r:prob}).
Для доказательства существования <<хорошего>> объекта подсчитывается, что <<плохих>> объектов  меньше, чем всех.   Кроме того, теоремы излагаются без технических усложнений, дающих незначительно более сильные оценки
(т.о. обычно теоремой Эрдеша и теоремой об оценке уклонения называются не теоремы \ref{t:erd} и \ref{t:rad}, а немного более сильные утверждения).
Также мы постарались хорошо структурировать доказательство теоремы \ref{t:rad} (т.е. разбить его на шаги, в частности, выделить красивую лемму \ref{l:one}).

Мы начинаем с олимпиадных задач, являющихся частными случаями указанных теорем.

\begin{problem}\label{p:ram} (a)  Среди любых 51 из 10 миллионов китежан имеется двое знакомых.
Обязательно ли найдется 51 китежанин, любые два из которых знакомы?

(b) Любые два из 1000 ученых переписываются по одной из четырех тем: географии, геологии, топографии и топологии.
Обязательно ли  найдутся 12 ученых, любые два из которых переписываются по одной и той же теме?

(c) Среди 1000 членов хурала выбрано несколько комиссий, в каждой из которых 3 человека.
Обязательно ли  найдется 10 членов хурала, из которых либо любые 3 образуют комиссию, либо любые 3 не образуют?
\end{problem}

Ответы --- нет.
Они являются частным случаем следующих теорем.

Назовем {\it $n$-кликой} ({\it $n$-антикликой}) полный (пустой) подграф на $n$ вершинах данного графа.
Для фиксированного $l$ назовем {\it $n$-гиперкликой} семейство всех $l$-элементых подмножеств данного подмножества из $n$ элементов.

\begin{theorem}[Эрдеш]\label{t:erd} (a) Для любого $n$ существует граф на $2^{\lfloor\frac{n-2}2\rfloor}$ вершинах, не имеющий ни $n$-клики, ни $n$-антиклики.

(b) Для любых $n,k$ существует раскраска ребер полного графа на $k^{\lfloor\frac{n-2}2\rfloor}$ вершинах в $k$ цветов, для которой нет одноцветной $n$-клики.

(c) Для любых $n,k,l$ существует раскраска всех $l$-элементных подмножеств множества из $k^{\lfloor(n-l+1)^{l-1}/l!\rfloor}$ элементов в $k$ цветов, в котором нет одноцветной $n$-гиперклики.
%$n$-элементного подмножества, все $l$-элементные подмножества которого одноцветны.
\end{theorem}

%При помощи индукции по $n\ge3$ построим набор комиссий в хурале из $2^{n-2}$ человек, для  которого среди
%любых $n$ человек найдется комиссия, а также найдется тройка человек, не образующих комиссию.
%Для $n=3$ построение очевидно.
%Пусть $X_n$ --- нужный хурал из $2^{n-2}$ человек.
%Обозначим через $Y_n$ копию этого хурала.
%Составим комиссию из тройки $A\subset X_n\sqcup Y_n$, если
%либо $A$ --- комиссия в $X_n$, либо $A$ --- комиссия в $Y_n$, либо $|A\cap X_n|=2$.
%Получим нужный набор комиссий в хурале  $X_n\sqcup Y_n$ из $2^{n-1}$ человек.

Следующая задача является частным случаем теоремы \ref{t:rad} об оценке уклонения для $n=1000$, $s=300$ и $a=150$ (поскольку $2^{150^2}>2^{10\cdot2000}>600^{2000}$).

\begin{problem}[\cite{PT}]\label{p:rad}
1000 пионеров города Новые Васюки вышли на парад. Известно, что пионеры ходят в 300 кружков.
Докажите, что Остап Бендер может раздать пионерам пилотки двух цветов (красные и синие) так, чтобы среди представителей одного кружка разность (по модулю) между количествами пионеров в красных пилотках и в синих пилотках не превосходила 150.
\end{problem}

\begin{theorem}[об оценке уклонения]\label{t:rad}
Пусть $M_1,\ldots,M_s$ --- семейство подмножеств множества $[n]:=\{1,\ldots,n\}$.
Если $2^{a^2}\ge(2s)^{2n}$, то существует раскраска множества $[n]$ в красный и синий цвета, для которой при любом $k\in[s]$ количества красных и синих элементов в $M_k$ отличаются менее чем на $a$.
%$\sqrt{2n\ln(2s)}$.
\end{theorem}

\begin{proof}[Доказательство ответов <<нет>> в задачах \ref{p:ram}.ac]
(a) Рассмотрим граф знакомств китежан.
Назовем граф с вершинами $1,2,\ldots,10^7$ {\it рамсеевским}, если в нем есть либо 51-клика, либо 51-антиклика.
Для доказательства существования нерамсеевского графа достаточно показать, что рамсеевских графов меньше,
чем всех графов с вершинами $1,2,\ldots,10^7$.
Количество последних равно $N:=2^{10^7(10^7-1)/2}$.
Клику на данной 51 вершине можно продолжить до $N/2^{51\cdot50/2}$ графов с вершинами $1,2,\ldots,10^7$.
Аналогичное справедливо для антиклики.
Количество 51-элементных подмножеств $10^7$-элементного множества равно ${10^7\choose51}<(10^7)^{51}=10^{357}$.
Поэтому количество рамсеевских графов меньше, чем
$$2\cdot10^{357}\frac{N}{2^{51\cdot25}}< \frac{N\cdot10^{357}}{2^{10\cdot5\cdot25}}< \frac{N\cdot10^{357}}{10^{3\cdot5\cdot25}}=N\cdot10^{357-375}<N.$$

(c) Назовем семейство трехчеловечных комиссий (среди 1000 членов хурала) {\it рамсеевским}, если в нем есть либо 10-гиперклика, либо 10-антигиперклика (10-антигиперкликой называется 10 человек, никакая тройка из которых не является комиссией).
Для доказательства существования нерамсеевского семейства трехчеловечных комиссий достаточно показать, что рамсеевских семейств меньше, чем всех семейств.
Количество последних равно $N:=2^{{1000\choose3}}$.
Гиперклику размера 10 данного подмножества из 10 членов хурала можно продолжить до $N/2^{{10\choose3}}$ семейств трехчеловечных комиссий.
Антигиперклику размера 10 данного подмножества из 10 членов хурала можно продолжить до $N/2^{{10\choose3}}$ семейств трехчеловечных комиссий.
Количество 10-элементных подмножеств 1000-элементного множества равно ${1000\choose10}<(1000)^{10}<2^{100}$.
Поэтому количество рамсеевских семейств трехчеловечных комиссий меньше, чем
$$2\cdot2^{100}\frac{N}{2^{{10\choose3}}} = \frac{N\cdot2^{101}}{2^{120}} < N.$$
\end{proof}

Ответ <<нет>> в задаче \ref{p:ram}.b и теорема \ref{t:erd} Эрдеша доказываются аналогично.
В частности, в начале доказательства теоремы \ref{t:erd}.a берем $r:=2^{\lfloor\frac{n-2}2\rfloor}$, тогда
$\binom{r}{n}<r^n\le 2^{(n^2-n-2)/2}$ при $n\ge2$.

\smallskip
{\it Перейдем к доказательству теоремы \ref{t:rad} об оценке уклонения} (и, тем самым, к решению задачи \ref{p:rad}).
Для подмножества $M$ множества $[n]$ и раскраски $x$ множества $[n]$ в два цвета обозначим через $\Delta_M(x)$ ({\it уклонение}) разность
%(без модуля)
между количествами элементов первого и второго цвета в $M$.

\begin{lemma}\label{l:one} Для любых $a>0$ и подмножества $M$ множества $[n]$ количество раскрасок $x$, для которых $\Delta_M(x)\ge a$, меньше $2^{n-\frac{a^2}{2n}}$.
%$2^ne^{-\frac{a^2}{2n}}$.
\end{lemma}

\begin{proof}[Доказательство теоремы \ref{t:rad}]
Применим лемму \ref{l:one} к любому $k\in[s]$.
Получим, что количество раскрасок $x$, для которых $\Delta_{M_k}(x)\ge a$, меньше $2^{n-\frac{a^2}{2n}}\le2^n/(2s)$.
В силу симметрии количество раскрасок $x$, для которых $\Delta_{M_k}(x)\le -a$, также меньше $2^n/(2s)$.
Поэтому количество ($k$-<<плохих>>) раскрасок $x$, для которых $|\Delta_{M_k}(x)|\ge a$, меньше $2^n/s$.
Тогда количество (<<плохих>>) раскрасок $x$, для которых найдется такое $k\in[s]$,
что $|\Delta_{M_k}(x)|\ge a$, меньше $s2^n/s=2^n$.
Следовательно, найдется <<хорошая>> раскраска.
\end{proof}

В доказательстве леммы \ref{l:one} используется следующий очевидный, но крайне полезный факт.

\begin{lemma}[Неравенство Маркова]\label{t:mar}
Для любых $w_1,\ldots,w_s>0$ количество тех $i$, для которых $w_i \ge 1$, не превосходит $w_1+\ldots+w_s$.

Более общо, для любых $a,w_1,\ldots,w_s>0$ количество тех $i$, для которых $w_i \ge a$, не превосходит
$\dfrac{w_1+\ldots+w_s}{a}$.
\end{lemma}

\begin{proof}[Доказательство леммы \ref{l:one}]
Ввиду неравенства Маркова (лемма \ref{t:mar}) для любого $\lambda>0$ имеем
$$
|\{x\ :\ \Delta_M(x)\ge a\}|\ =\ |\{x\ :\ \lambda\Delta_M(x)\ge\lambda a\}|\ =
\ |\{x\ :\ e^{\lambda\Delta_M(x)}\ge e^{\lambda a}\}|\ \le\ e^{-\lambda a} \sum_x e^{\lambda\Delta_M(x)}.$$
Раскраску подмножества $A\subset [n]$ в два цвета будем считать отображением $A\to\{-1,1\}$.
Тогда
$$
\sum_{x\in\{-1,1\}^{[n]}} e^{\lambda\Delta_M(x)} =
\sum_{x\in\{-1,1\}^{[n]}} \prod_{j \in M} {e^{\lambda x(j)}} =
2^{n-|M|}\sum_{y\in\{-1,1\}^M} \prod_{j \in M} {e^{\lambda y(j)}} =
2^{n-|M|}(e^{\lambda}+e^{-\lambda})^{|M|}.
$$
Здесь

$\bullet$ первое равенство верно, поскольку $\Delta_M(x)=\sum_{j\in M}x(j)$;

$\bullet$ второе равенство верно, поскольку каждую раскраску подмножества $M$ в два цвета  можно продолжить до $2^{n-|M|}$ раскрасок множества $[n]$.

$\bullet$ третье равенство верно, поскольку при раскрытии скобок в правой части получается левая часть.

Поэтому
$$|\{x\ :\ \Delta_M(x)\ge a\}| \le
e^{-\lambda a}2^n\left(\frac{e^\lambda + e^{-\lambda}}{2}\right)^{|M|} \overset{(*)}<
2^ne^{|M|\frac{\lambda^2}{2} - \lambda a} \leq
2^ne^{n\frac{\lambda^2}{2} - \lambda a} \overset{(**)} = 2^ne^{-\frac{a^2}{2n}} < 2^{n-\frac{a^2}{2n}}.
$$
Здесь

$\bullet$ неравенство $(*)$ верно, поскольку
%Покажем, что $e^\lambda+e^{-\lambda} < 2e^{\frac{\lambda^2}2}$ .
$$
e^{\lambda}+e^{-\lambda} =
\sum_{k=0}^{\infty}\dfrac{\lambda^k}{k!}+\sum_{k=0}^{\infty}\dfrac{(-1)^k\lambda^k}{k!} =
2\sum_{k=0}^{\infty}\dfrac{\lambda^{2k}}{(2k)!} <
2\sum_{k=0}^{\infty}\dfrac{\lambda^{2k}}{2^kk!} = 2e^{\frac{\lambda^2}2}.
$$

$\bullet$ равенство $(**)$ получается подстановкой $\lambda=a/n$.
\end{proof}

\begin{remark}\label{r:prob} Часто для доказательства существования <<хорошего>> объекта {\it подсчитывается}, что <<плохих>> объектов  меньше, чем всех.
Такие доказательства существования можно излагать на вероятностном языке (т.е. построив дискретное вероятностное пространство с равновероятными элементарными событиями).
К сожалению, такое изложение недоступно для многих студентов.
Ибо многие не в состоянии сформулировать необходимое для доказательства определение (дискретного)  вероятностного пространства.
А многим не хватает математической культуры даже для осознания того, что такого пространства нет в формулировке теоремы, поэтому его построение --- часть доказательства (<<дополнительное построение>>).

Заблуждение о наличии <<явных недостатков>> простых изложений по сравнению с наукообразными настолько распространено, что полезно сделать следующие замечания.
Хотя подсчет <<плохих>> объектов позволяет обойтись без вероятностного языка, сущностные
(т.е. не связанные с языком изложения) детали доказательства при обоих изложениях одни и те же.
Все простые соображения, которые используются выше, нужны и при вероятностном изложении. 
Использование терминов (например, определения матожидания) не упрощает и не сокращает выкладки, а лишь позволяет говорить <<по [очевидному] свойству матожидания>> вместо <<очевидно>>.
Для большинства читателей возможность самостоятельно проделать выкладки возрастает с простотой языка, подготавливающего эти выкладки.

%\footnote{Для читателя, знакомого с вероятностным изложением, заметим, что доказательство леммы \ref{l:one}
%на вероятностном языке включает в себя доказательство мультипликативности  матожидания.
%Для последнего может оказаться необходимым переход между римановским и лебеговским определениями матожидания --- %если подсчет матожидания проведен не по тому определению, по которому удобно доказывать мультипликативность.}
%В доказательстве леммы 5 есть цепочка равенств с суммами и произведениями, последнее из них 
%(которое получается раскрытием скобок) как раз соответствует мультипликативности. 
%Можно сказать, что действительно мультипликативности след простыл

Тем не менее, развивать вероятностную интуицию крайне полезно.
Применение дискретных вероятностных пространств с элементарными событиями, уже не являющимися равновероятными, --- мощный метод комбинаторики, \emph{вероятностный метод} \cite{AS, R3}.
В качестве пропедевтики вероятностного метода полезны доказательства существования, использующие подсчет
(подобно вышеприведенным).
После их изучения разумно потренироваться излагать эти доказательства на вероятностном языке
\cite[решения к п. 1.6 <<подсчет двумя способами>>]{GDI}.
(Точно так же, как в начале освоения теории Галуа разумно потренироваться излагать на ее языке уже известные решения квадратного уравнения и уравнений 3-й и 4-й степени.)
Другие примеры естественного выращивания вероятностного языка приведены, например, в \cite[\S22]{ZSS}, \cite{IRS}, \cite[\S6.2]{GDI}.
\end{remark}

%Поэтому если и найдутся знакомые с такими книгами и лекциями, то они не будут иметь преимущества перед теми,
%кто способен построить простой пример или сообразить, что для доказательства существования хорошего
%объекта достаточно доказать, что плохих меньше, чем всех.

\end{document}